\documentclass{amsproc}

\newtheorem{Theorem}{Theorem}

\newtheorem{theorem}{Theorem}[section]
\newtheorem{lemma}[theorem]{Lemma}
\newtheorem{proposition}[theorem]{Proposition}

\theoremstyle{definition}
\newtheorem{definition}[theorem]{Definition}

\theoremstyle{remark}

\numberwithin{equation}{section}

\newcommand{\abs}[1]{\lvert#1\rvert}

\newcommand{\R}{{\mathbb R}}

\newcommand{\FF}{{\mathcal F}}

\newcommand{\Int}{{{\rm Int}\,}}

\newcommand{\LL}{{\mathcal L}}

\newcommand{\BB}{{\mathcal B}}

\begin{document}
\title[Continuous leafwise harmonic functions]
{Conitinuous leafwise harmonic functions on codimension one
transversely isometric foliations}


\author{Shigenori Matsumoto}
\address{Department of Mathematics, College of
Science and Technology, Nihon University, 1-8-14 Kanda, Surugadai,
Chiyoda-ku, Tokyo, 101-8308 Japan
}
%\curraddr{Department of Mathematics, College of
%Science and Technology, Nihon University, 1-8-14 Kanda, Surugadai,
%Chiyoda-ku, Tokyo, 101-8308 Japan}
\email{matsumo@math.cst.nihon-u.ac.jp
}
\thanks{The author is partially supported by Grant-in-Aid for
Scientific Research (C) No.\ 25400096.}
%    General info


\date{\today }
\begin{abstract}
Let $\FF$ be a codimension one foliation on a closed manifold $M$
which admits a transverse dimension one Riemannian foliation.
Then any continuous
leafwise harmonic functions are shown to be constant on leaves.
\end{abstract}

\maketitle

\section{Introduction}

Let $M$ be a closed $C^2$ manifold, and let $\FF$ be a continuous
leafwise $C^2$ foliation on $M$. This means that $M$ is covered by
a finite union of continuous foliation charts and the transition
functions
are continuous, together with their leafwise partial derivatives up to
order 2. Let $g$ be a continuous leafwise\footnote{This means
that the leafwise
partial derivatives up to order 2
of the components of $g$ in each foliation chart is continuous in the chart.}
$C^2$ leafwise Riemannian metric.
In this  paper such a triplet $(M,\FF,g)$ is simply refered to as
a {\em leafwise $C^2$ foliations.} For simplicity, we assume
throughout the paper that the manifold
$M$ and the foliation $\FF$ are oriented.
For a continuous leafwise\footnote{The leafwise
partial derivatives of $h$ up to order 2
in each foliation chart is continuous in the chart.}
$C^2$ real valued function $h$ on $M$, the leafwise Laplacian $\Delta h$ is defined
using the leafwise metric $g$. It  is a continuous function of $M$.

\begin{definition}
A continuous leafwise $C^2$ function $h$ is called {\em leafwise
 harmonic}
if $\Delta h=0$.
\end{definition}

\begin{definition}
A leafwise $C^2$ foliation
 $(M,\FF,g)$ is called {\em Liouville} if any continuous
leafwise harmonic function is leafwise constant.
\end{definition}

As an example, if $\FF$ is a foliation by compact leaves, then
$(M,\FF,g)$ is Liouville. Moreover there is an easy observation.

\begin{proposition}\label{minimal}
If $\FF$ admits a unique minimal set, then $(M,\FF,g)$ is Liouville.
\end{proposition}

This can be seen as follows. Let $m_1$ (resp.\ $m_2$) be the maximum (resp.\
minimum)
value of the continuous leafwise harmonic function $h$ on $M$. Assume
$h$ takes the maximum value  $m_1$ at $x\in M$. Then by the maximum
principle, $h=m_1$ on the leaf $F_x$ which passes through $x$.
Now the closure of $F_x$ contains the unique minimal set $X$.
Therefore $h=m_1$ on $X$. The same argument shows that $h=m_2$ on $X$,
finishing the proof that $h$  is constant on $M$.

A first example of non-Liouville foliations is obtained by R. Feres and
A. Zeghib in a beautiful and simple construction \cite{FZ}.
It is a foliated $S^2$-bundle over a hyperbolic surface, with two
compact leaves. There are also examples in codimension one.
B. Deroin and V. Kleptsyn \cite {DK} have shown that
a codimension one foliation $\FF$ is non-Louville if $\FF$ is
transversely $C^1$, admits no transverse invariant measure and
possesses more than one minimal sets, and they have constructed such
a foliation.

A codimension one foliation $\FF$ is called {\em $\R$-covered} if the leaf
space
of its 
lift to the universal covering space is homeomorphic to $\R$.
See \cite{F} or \cite{FFP}.
It is shown in \cite{F} and \cite{DKNP} that an
$\R$-covered
foliation without compact leaves admits a unique minimal set.
Therefore the above example of  a codimension one non-Liouville foliation
is not $\R$-covered.
This led the authors of \cite{FFP} to the study of Liouville property
for $\R$-covered foliations. The purpose of this paper is to generalize
a result of \cite{FFP}.

\begin{definition}
A codimension one leafwise $C^2$ foliation $(M,\FF,g)$ is called 
{\em transversely isometric} if there is a continuous dimension one
foliation
$\phi$ transverse to $\FF$ such that the holonomy map of $\phi$
sending a (part of a) leaf of $\FF$ to another leaf is $C^2$ and preserves
the leafwise metric $g$.
\end{definition}

Notice that a transversely isometric foliation is $\R$-covered.
Our main result is the following.

\begin{Theorem} \label{T}
A leafwise $C^2$ transversely isometric
codimension one foliation is Liouville.
\end{Theorem}

In \cite{FFP}, the above theorem is proved in the case where
the leafwise Riemannian metric is negatively curved. 
Undoubtedly this is the most important case. But 
the general case may equally be of interest.

If a transversely isometric foliation $\FF$ does not admit a compact leaf,
then, being $\R$-covered,
it admits a unique minimal set, and Theorem \ref{T} holds true
by Proposition \ref{minimal}.
Therefore we only consider the case where $\FF$ admits a compact leaf.
In this case the union $X$ of compact leaves is closed. Let
$U$ be a connected component of $M\setminus X$, and let $N$ be
the metric completion of $U$. Then $N$ is a foliated interval bundle,
since the one dimensional transverse foliation $\phi$ is Riemannian.

Therefore we are led to consider the following situation.
Let $K$ be a closed $C^2$ manifold of dimension $\geq 2$, equipped with
a $C^2$ Riemannian metric $g_K$. Let $N=K\times I$, where $I$ is
the interval $[0,1]$. Denote by $\pi:N\to K$ the canonical projection.
Consider a continuous foliation $\LL$ which
is transverse to the fibers $\pi^{-1}(y)$, $\forall y\in K$.
Although $\LL$ is only continuous, its leaf has a $C^2$ differentiable
structure as a covering space of $K$ by the restriction of
$\pi$. Also $\LL$ admits
a leafwise Riemannian metric $g$ obtained as the lift of $g_K$ to each
leaf by $\pi$. 
Such a triplet $(N,\LL,g)$ is called a {\em leafwise $C^2$ foliated $I$-bundle}
in this paper.
Now Theorem \ref{T} reduces to the following theorem.

\begin{Theorem} \label{main}
Assume a leafwise $C^2$ foliated $I$-bundle 
$(N,\LL,g)$ does not admit a compact leaf in the interior $\Int(N)$.
Then any continuous leafwise harmonic function is constant on $N$.
\end{Theorem}

An analogous result for random discrete group actions on the interval
was obtained in \cite{FR}.

The rest of the paper is devoted
to the proof of Theorem \ref{main}. The proof is by absurdity.
Throughout the paper, $(N,\LL,g)$
denotes a leafwise $C^2$ foliated $I$-bundle without interior compact leaves,
and we assume 
that there is a
continuous leafwise harmonic function $f$ such that $f(K\times\{i\})=i$,
$i=0,1$. As is remarked in \cite{FFP}, this is not a loss of generality.
Also notice that for any point $x\in\Int(N)$, we have $0<f(x)<1$.

\section{Preliminaries}
In this section, we recall fundamental facts about Brownian motions,
needed in the next section.

Let us denote by $\Omega$ the space of continuous leafwise paths
$\omega:[0,\infty)\to N$. For any $t\geq0$, a random variable
$X_t:\Omega\to N$ is defined by $X_t(\omega)=\omega(t)$.
Let $\BB$ be the $\sigma$-algebra of $\Omega$
generated by $X_t$ ($0\leq t<\infty$).
 As is well known, easy to show, $\BB$ coincides with the
$\sigma$-algebra
generated by the compact open topology on $\Omega$.
A bounded function $\phi:\Omega\to \R$ is called a {\em Borel function}
if $\phi$ is $\BB$-measurable, i.\ e.\ if
for any Borel subset $B\subset\R$, the inverse image $\phi^{-1}(B)$
belongs to $\BB$.

For any point $x\in N$, the Wiener probability measure $P^x$ is defined
using the leafwise Riemannian metric $g$. Notice that
$P^x\{X_0=x\}=1$.
For any bounded Borel function $\phi:\Omega\to\R$,
the expectation of $\phi$ w.\ r.\ t.\ $P^x$ is denoted by $E^x[\phi]$.
The following proposition is well known.

\begin{proposition}\label{p0}
Let $f$ be a bounded Borel function defined on a leaf
$L$ of $\LL$ Then $f$ is harmonic on $L$ if and only if
\begin{equation}\label{e1}
E^x[g(X^t)]=g(x),\ \ \forall t\geq 0, \ \forall x\in L.
\end{equation}
\end{proposition}

Let $\FF$ be the completion of $\BB$ by the measure $P^x$.
For any $t\geq0$,
let $\BB_t$ be the $\sigma$-algebra generated by
$X_s$ ($0\leq s\leq t$). Its completion is denoted by $\FF_t$.
Notice that unlike $\BB$ and $\BB_t$, $\FF$ and $\FF_t$ depend
strongly on $x\in M$, but we depress the
dependence in the notation.

The following fact is well known and can be shown easily using
the Radon-Nikod\'ym theorem. 

\begin{proposition}\label{p1}
Given any bounded $\FF$-measurable function $\phi:\Omega\to\R$
and $t>0$,
there is a unique bounded $\FF_t$-measurable function $\phi_t$ such that
for any bounded $\FF_t$-measurable function $\psi$, we have
$$
E^x[\psi\phi_t]=E^x[\psi\phi].$$
\end{proposition}

The $\FF_t$-measurable function $\phi_t$ is called the {\em conditional
expectation of $\phi$ with respect to $\FF_t$} and is denoted by
$E^x[\phi\mid\FF_t]$.

For any $t>0$, let $\theta_t:\Omega\to\Omega$ be the shift map by $t$,
defined by
$$
\theta_t(\omega)(s)=\omega(t+s), \forall s\in[0,\infty).$$

The following proposition is known as the Markov property.
See for example \cite{O}.

\begin{proposition} \label{p2}
Let $\phi:\Omega\to\R$ be a bounded $\FF$-measurable function. Then we have
\begin{equation}\label{e2}
E^x[\phi\circ\theta_t\mid\FF_t]=E^{X_t}[\phi].
\end{equation}
\end{proposition}

A family $\phi_t$ of uniformly bounded $\FF_t$-measurable functions, $t\geq0$, is called
a {\em bounded $P^x$-martingale} if 
\begin{equation} \label{e3}
E^x[\phi_{t+h}\mid\FF_t]=\phi_t,\ \ \  \forall t\geq 0, \ \forall h>0.
\end{equation}

We have the following martingale convergence theorem.
See \cite{O}, Appendix C.
\begin{proposition}\label{p3}
Let $\{\phi_t\}$ be a bounded $P^x$-martingale. Then there is a bounded
 $\FF$-measurable
function $\phi$ such that $\phi_t\to\phi$ as $t\to\infty$,
$P^x$-almost surely. 
\end{proposition}
 
We shall raise two applications of the above facts, which will be useful 
in the next section.
Let $f:N\to[0,1]$ be the continuous leafwise harmonic function defined at the
end of Section 1.
Then $f(X_t)$ is a Borel function defined on $\Omega$.

\begin{lemma}\label{l1}
For any $x\in N$, there is an $\FF$-measurable function
$\phi:\Omega\to[0,1]$ such that $f(X_t)\to\phi$ as $t\to\infty$, $P^x$-almost surely.
\end{lemma}

{\sc Proof.} We only need to show that $f(X_t)$ is a bounded
$P^x$-martingale. For this, we have
$$
E^x[f(X_{t+h})\mid\FF_t]=E^x[f(X_{h})\circ\theta_t\mid\FF_t]=E^{X_t}[f(X_h)]=f(X_t),$$
where the second equality is by the Markov propery, and the last
by the leafwise harmonicity of $f$.
\qed

\begin{lemma}\label{l2}
Let $\phi:\Omega\to[0,1]$ be a Borel function such that
 $\phi\circ\theta_t=\phi$ for any $t>0$. Then
$x\mapsto E^x[\phi]$ is a harmonic function on each leaf $L$ of $\LL$.
\end{lemma}

{\sc Proof.} By Proposition \ref{p0}, we only need to show that
$$
E^x[E^{X_t}[\phi]]=E^x[\phi],\ \ \forall x\in L,\ \forall t>0.$$
But we have
$$E^x[E^{X_t}[\phi]]=E^x[E^x[\phi\circ\theta_t\mid\FF_t]]=
E^x[E^x[\phi\mid\FF_t]]=E^x[\phi].
$$
\qed

\section{Proof of Theorem \ref{main}}
 
Again let $f$ be a continuous leafwise harmonic function defined
at the end of Section 1.
A probability measure $\mu$ on $N$ is called {\em stationary} 
if $\langle\mu,\Delta h\rangle=0$ for any continuous leafwise $C^2$
function $h$.

\begin{proposition}\label{p31}
There does not exist a stationary measure $\mu$ such that
 $$\mu(\Int(N))>0.$$
\end{proposition}

{\sc Proof.} Denote by $X$ the union of leaves on which
$f$ is constant. The subset $X$ is closed in $N$.
L. Garnett \cite{G} has shown that  $\mu(X)=1$ for any stationary
measure
$\mu$.
Therefore if $\mu(\Int(N))>0$, there is a leaf $L$ in $\Int(N)$ on which
$f$ is constant. But since we are assuming that there is no interior
compact leaves, the closure of $L$ must contain  both boundary
components
of $N$. A contradiction to the continuity of $f$.
\qed

\bigskip
Given $0<\alpha<1$, let $V=K\times (\alpha,1]$, a neighbourhood of the
upper boundary component $K\times\{1\}$. Let
$\Omega_V$ be the subset of $\Omega$ defined by
$$
\Omega_V=\{X_{t_i}\in V, \ \exists t_i\to\infty\},$$
and let $\phi$ be the characteristic function of $\Omega_V$. Clearly
$\phi$ is a Borel function on $\Omega$ and satisfies
$\phi\circ\theta_t=\phi$
for any $t>0$. Thus by Lemma \ref{l2}, the function $p:M\to[0,1]$
defined by $p(x)=E^x[\phi]$
is leafwise harmonic. 

Another important feature of the function $p$ is that $p$ is nondecreasing along the
fiber $\pi^{-1}(y)$,
$\forall y\in K$, since our leafwise Brownian motion is synchronized,
i.\ e,\ it is the lift of
the Brownian motion on $K$. Notice that $p=i$ on $K\times\{i\}$,
$i=0,1$.

\begin{lemma}\label{l31}
The function $p$ is constant on $\Int(N)$.
\end{lemma}

The proof is the same as the proof of Proposition 9.1 of \cite{FFP}.
In short, if we assume the contrary, 
we can construct a stationary measure $\mu$ such
that $\mu(\Int(N))>0$, which is contrary to Proposition \ref{p31}.
The proof is included in Section 4 for completeness.

\begin{lemma}\label{l32}
The function $p$ is 1 on $\Int(N)$.
\end{lemma}

{\sc Proof.} Assume $p<1$ on $\Int(N)$. 
For any $x\in N$, $P^x$-almost surely the limit 
$\displaystyle \phi=\lim_{t\to\infty}f(X_t)$ exists by Lemma \ref{l1}.
Choose a constant $0<a<1$ so that $f^{-1}[a,1]$ is contained in $V$.
Then we have for any $x\in\Int(N)$,
$$
q_x=P^x\{\phi\geq a\}\leq p.$$
Therefore
\begin{eqnarray*}
E^x[\phi]&\leq &1\cdot P^x\{\phi\geq a\}+a\cdot P^x\{\phi< a\}\\
&=&q_x+a(1-q_x)=a+(1-a)q_x\leq a+(1-a)p<1.
\end{eqnarray*}

Now by the dominated convergence theorem, we have
$$
E^x[\phi]=E^x[\lim_{t\to\infty}f(X_t)]=\lim_{t\to\infty}E^x[f(X_t)]
=\lim_{t\to\infty} f(x)=f(x).$$
Since $x$ is an arbitrary point in $\Int(N)$, this shows that
$f$ cannot take value greater than $a+(1-a)p$ in $\Int(N)$,
contradicting the continuity of $f$.
\qed

\bigskip
{\sc Proof of Theorem \ref{main}.} Fix $x\in\Int(N)$.
Now for any small neighbourhood $V$ of $K\times\{1\}$, we have
$P^x(\Omega_V)=1$. This shows that 
$\displaystyle \limsup_{t\to\infty}f(X_t)=1$, $P^x$-almost surely.
Likewise considering small neighbourhoods of $K\times\{0\}$,
we have
$\displaystyle \liminf_{t\to\infty}f(X_t)=0$, $P^x$-almost surely.
But this contradicts Lemma \ref{l1}. We are done with the proof of
Theorem \ref{main}.

\section{Proof of Lemma \ref{l31}}

The projection $\pi:N\to K$ has a distinguished role since the leafwise
metric $g$ is the lift of $g_K$ by $\pi$.
Let $B$ be an open ball in the base manifold $K$, and consider an
inclusion  $\iota:B\times I\subset M$ such that
$\iota(\{y\}\times I)=\pi^{-1}(y)$, $\forall y\in B$ and 
$\iota(B\times\{t\})$ is contained
in a leaf of $\LL$, $\forall t\in I$. 
Such an inclusion $\iota$, or its domain $B\times I$, is called a {\em
distinguished chart}.

Using a distinguished chart, let us define $\hat p:N\to [0,1] $ 
by
\begin{eqnarray*}
\hat p(y,t)&=&\lim_{h\downarrow0}p(y,t+h)\ \mbox{ if } t<1,\\
\hat p(y,1)&=&p(y,1).
\end{eqnarray*}

Of course $\hat p$ does not depend on the choice of the distinguished charts,
and is defined on the whole $N$. 

The function $\hat p$ is right semicontinuous on each fiber
$\pi^{-1}(y)$,
$\forall y\in K$, and is leafwise harmonic just as $p$.
Since $\hat p$ is nondecreasing
on
each fiber $\pi^{-1}(y)$, $y\in K$, we get a probability measure $\nu_y$
on $\pi^{-1}(y)$ by
\begin{eqnarray*}
\nu_y(\{y\}\times(s,t])&=&\hat p(y,t)-\hat p(y,s),\\
\nu_y(\{y\}\times\{0\})&=&\hat p(y,0)-p(y,0).
\end{eqnarray*}
for any $0\leq s<t\leq 1$.
Again $\nu_y$ does not depend on the choice of the distinguished charts.

Our goal is to show that $\nu_y(\Int(\pi^{-1}(y)))=0$ for any $y\in K$.
Clearly this shows that $\hat p$, and hence $p$, is constant on $\Int(N)$.
Since we do not use the definition of $p$ in what follows, 
we write $p$ for $\hat p$ for simplicity. 

\bigskip
We first show that $\nu_y$ does not have an atom in $\Int(\pi^{-1}(y))$,
$\forall y\in K$. Assume on the contrary that
there is $y_0\in K$ such that $\nu_{y_0}$ has an atom in
$\Int(\pi^{-1}(y_0))$.
Then in a distinguished chart, there is $t_0\in(0,1)$
such that
$$
q(y_0,t_0)=p(y_0,t_0)-\lim_{t\uparrow t_0}p(y_0,t)>0.$$
Define a positive function $q$ on the plague $B\times\{t_0\}$
by
$$
q(y,t_0)=p(y,t_0)-\lim_{t\uparrow t_0}p(y,t).$$
Using other distinguished charts,
we can define $q$  on the whole leaf $L$
which passes through $(y_0,t_0)$. 
The function $q$ is positive harmonic on $L$.
Clearly we have
$$
\sum_{x\in L\cap\pi^{-1}(y)}q(x)\leq 1,$$ 
for any $y\in K$.

Define a measure $\mu$ on $N$ by
$$\mu=\int_K(\sum_{x\in L\cap\pi^{-1}(y)}q(x)\delta_x)dy,
$$
where $\delta_x$ is
the Dirac mass at $x$, and $dy$ denotes the volume form on $K$. 
Then according to a criterion in \cite{G},
$\mu$ is a stationary measure, contradicting Proposition \ref{p31}.

\bigskip
Now $p$ is continuous on $\Int(N)$, since there is
no atom of $\nu_y$.
Next assume that there is $y_0\in K$ and an interval $J$ in
$\pi^{-1}(y_0)$
such that $\nu_{y_0}(J)=0$.
Let $J$ be a maximal such interval and write it
as $J=\{y_0\}\times[a,b]$ in a distinguished chart
$B\times I$ such that $y_0\in B$. Then in that chart the function 
$y\mapsto p(y,b)-p(y,a)$ is nonnegative harmonic, and takes value $0$ 
at $y_0$. Therefore we have 
$$
p(y,b)=p(y,a),\ \forall y\in B.$$
This equality holds also in neighbouring distinguished charts and therefore
all over $N$.
Thus if we delete the saturation of $\Int(J)$ from $N$, the function $p$ is
still well defined and continuous.
Deleting all such open saturated sets,
we get a  new manifold, still denoted by $N$, and a new
foliation, still denoted by $\LL$.

\bigskip
For the new foliation 
$\LL$, the function
$p$ is continuous and strictly increasing on the interior of each
fiber $\pi^{-1}(y)$. It may not be continuous on the boundary.
However it is possible to
 extend the function $p\vert_{\Int(N)}$ to the boundary 
by continuity, thanks to the leafwise harmonicity
of $p$.
The new function $p$ is still constant on each boundary
component of $N$. After normalization, 
one may assume that $p=i$ on $K\times\{i\}$, $i=0,1$.
Let $B$ be an open ball in $K$ centered at $y_0$. Define a special
kind of distinguished charts
$B\times I$ by using the value $p$ along the fiber $\pi^{-1}(y_0)$.
That is, $p(y_0,t)=t$, $\forall t\in I$. 
Notice that $\nu_{y_0}$ is the Lebesgue measure $dt$ in this chart.

Given two values 
$0\leq t'<t''\leq 1$, the function $q(y)=p(y,t'')-p(y,t')$ is
a positive harmonic function on $B$. 
By the Harnack inequality, there is $C_1>0$
independent of $t'$ and $t''$ such that
$$
\abs{\log q(y)-\log q(y_0)}<C_1d(y,y_0)<C_2.
$$
Since $q(y_0)=t''-t'$, this shows that there is $C>1$ such that
\begin{equation}\label{e51}
C^{-1}\abs{t''-t'}\leq \abs{p(y,t'')-p(y,t')}\leq C\abs{t''-t'},
\end{equation}
for any $y\in B$ and $0\leq t'<t''\leq 1$.
Then $t\mapsto p(y,t)$ is Lipschitz, and thus differentiable $dt$-almost
everywhere. Precisely we have
\begin{eqnarray}\label{e52}
& &\mbox{for any $y\in B$, there is a Lebesgue full measure set $I_y$ of $I$
such that}
\\
& &\mbox{the partial derivative $p_t(y,t)$ exists for any $t\in I_y$.}
\nonumber
\end{eqnarray}

On the other hand, since 
the space of harmonic functions taking values in $[C^{-1},C]$
is compact,
(\ref{e51})
shows that the upper partial derivative $p_t^+(y,t)$ defined by
$\limsup$
and the lower
partial derivative $p_t^-(y,t)$
exist for any $(y,t)\in B\times I$, and is a harmonic function of $y$.
Notice also that $$C^{-1}\leq p_t^{\pm}(y,t)\leq C.$$

By (\ref{e52}), we have $p_t^+(y,t)=p^-_t(y,t)$ for any $t\in I_y$.
Then by Fubini, there is a Lebesgue full measure set $I_*$ of $I$
such that for any $t\in I_*$, 
\begin{equation}\label{e53}
p_t^+(y,t)=p_t^-(y,t) 
\end{equation}
holds for $dy$-almost all $y\in B$. But then (\ref{e53}) holds for
any $y\in B$, since $p^\pm_t(y,t)$ is harmonic in $y$. 

Writing the common value by $p_t(y,t)$ as usual, 
$\forall t\in I_*$, we
can define a measure $\mu$ on $B\times I$ by
$$
\mu=\int_{I_*}(p_t(y,t)dy)dt.
$$
This measure does not depend on the choice of the center $y_0$ of $B$.
Thus $\mu$ can be defined on the whole $N$.
Again by a criterion in \cite{G}, $\mu$ is a stationary measure. This contradicts Proposition
\ref{p31}, finishing the proof of Lemma \ref{l31}.

\end{document}